\documentclass[12pt]{amsart}
\newcommand{\g}{\mathfrak}
\newcommand{\cal}{\mathcal}

\newcommand{\gtg}{\g g}

\newcommand{\hgtg}{\hat{\gtg}}
\newcommand{\tgtg}{\tilde{\gtg}}

\newcommand{\gtsl}{{\g sl}}

\newcommand{\gtn}{{\g n}}
\newcommand{\gtnp}{\gtn_{+}}
\newcommand{\gtnm}{\gtn_{-}}
\newcommand{\gth}{{\g h}}
\newcommand{\nc}{{\bf C}}
\newcommand{\nq}{{\bf Q}}
\newcommand{\tens}{\dot{\otimes}}
\newcommand{\hg}{\hat{\Gamma}}

\newcommand{\nz} {{\bf Z}}

\newcommand{\cp}{{\bf CP^{1}}}

\newcommand{\co}{{\cal O}}
\newcommand{\tco}{{\tilde{\co_{k}}}}

\newcommand{\ca}{{\cal A}}

\newtheorem{theorem}{Theorem}[subsection]
\newtheorem{proposition}[theorem]{Proposition}

\newtheorem{remark}[theorem]{Remark}

\newtheorem{lemma}[theorem]{Lemma}
\newtheorem{corollary}[theorem]{Corollary}

\title[annihilating ideals and tilting functors]
{annihilating ideals and tilting functors}
\author{Igor B. Frenkel}
\address{Department of Mathematics, Yale University}
\author{Feodor Malikov}
\address{Department of Mathematics, University of Southern
California}

\normalsize 
\divide\textheight by \baselineskip
\multiply\textheight by \baselineskip
\begin{document}

\begin{abstract}

We use Kazhdan-Lusztig tensoring to, first, describe annihilating ideals
of highest weight modules over an affine Lie algebra in terms of the
corresponding VOA and, second, to classify tilting functors, an affine
analogue of projective functors known in the case of a simple Lie algebra.

\end{abstract}

\maketitle

 \section{introduction}

This paper grew out of an attempt to carry   the classical theory of
Harish-Chandra bimodules over to the case of an affine Lie algebra.
As we shall explain below, the very definition of a Harish-Chandra bimodule over
 an affine Lie
algebra is not obvious. We do not propose such a  definition, but we affinize
the notion of {\em projective functor}, an important tool of the classical
theory. In this way we get a semi-simple monoidal category, whose classical counterpart
is the subcategory of projective Harish-Chandra bimodules. The Grothendieck ring of this
category is  shown to be isomorphic to the group algebra of the corresponding affine
 Weyl group.

It is well-known that when studying representation theory of affine Lie algebras it is
important to distinguish between the case when the central charge is positive and the case when
it is negative. The results we have just discussed are all valid when the latter is the case.
Another part of the classical theory revolving
 about Duflo's theorem on primitive ideals
of the universal enveloping algebra generalizes to
 affine Lie algebras when the former is the
case. The generalization is non-trivial: we establish a bijection between the ideals
of the vertex operator
algebra   attached to a given affine Lie algebra and submodules of a Weyl module with a regular
dominant highest weight. This is an analogue of
 the Dixmier conjecture proved by Joseph \cite{joseph} which classifies
ideals of the universal enveloping algebra in terms of submodules of a
 Verma module with a regular dominant highest weight.

The main idea employed in this paper is to replace the functor of tensor product with a
finite-dimensional module with the functor of Kazhdan-Lusztig tensoring with a Weyl module,
or with a module having a filtration by Weyl modules. To make things clearer,
 we first review a functorial approach to 
Harish-Chandra
bimodules over a simple Lie algebra following the beautiful paper by
Bernstein and S.Gelfand \cite{bernst_gelf}.

\subsection{Representations of complex groups. }
\label{reviewofbernstgelf}
 Let $\ca$ be a category.
 Then
one can consider the category $Funct(\ca)$ of functors on $\ca$, 
objects
being functors, morphisms being natural transformations of functors.
In general,  there is
no reason to think that $Funct(\ca)$ is abelian even if
$\ca$ is so. Here is, however, an important example when $Funct(\ca)$ contains
an abelian complete subcategory.

Let $Mod(\gtg)$ be the category of  modules over a simple complex Lie algebra
$\gtg$ and $Mod(\gtg-\gtg)$ the category of $\gtg$-bimodules. 
``Module'' will always mean a space carrying a left action of $\gtg$; 
``bimodule'' will always mean a space carrying a left and a right action
commuting with each other.
Any $H\in  Mod
(\gtg-\gtg)$ gives rise to the functor
\[\Phi_{H}: Mod(\gtg)\rightarrow Mod(\gtg);\; \Phi_{H}(M)=H\otimes_{\gtg}M.\]

It is well-known that 
\[Hom_{Mod(\gtg-\gtg)}(H_{1},H_{2})=Hom_{Funct(Mod(\gtg))}(\Phi_{H_{1}},
\Phi_{H_{2}}).\]
\begin{sloppypar}
Therefore $Mod(\gtg-\gtg)$ is a complete abelian subcategory of $Funct(Mod(\gtg))$.
\end{sloppypar}

Any $\gtg$-bimodule is a $\gtg$-module with respect to the diagonal action
(that is, the left action minus the right action). A Harish-Chandra
 bimodule is
a finitely generated bimodule
 such that under the diagonal action it decomposes in a direct
sum of finite dimensional $\gtg$-modules occurring with finite multiplicities.
Consider the category of Harish-Chandra bimodules $HCh$, and  $\co_{\gtg}$,
the Bernstein-Gelfand-Gelfand category
of  $\gtg$-modules. The condition imposed
on the diagonal action ensures that if $H$ is a Harish-Chandra
bimodule, then $\Phi_{H}$ preserves $\co_{\gtg}$.
Therefore the construction we just discussed gives an
embedding $HCh\hookrightarrow Funct(\co_{\gtg})$ as a complete subcategory.

Further, indecomposable 
projective Harish-Chandra bimodules are exactly those corresponding
to direct indecomposable summands of the functor of tensoring by
a finite dimensional $\gtg$-module $V$:
\begin{equation}
\label{projfunctdef}
V\otimes ?: \co_{\gtg}\rightarrow\co_{\gtg},\; M\mapsto V\otimes M.
\end{equation}
Such functors are naturally called {\em projective}. To review their classification    obtained in
\cite{bernst_gelf}, we need some facts of representation theory.  

     $\co_{\gtg}$ admits the direct product decomposition with
respect to the
action of the center of the universal enveloping $U(\gtg)$:
  \[ \co_{\gtg}=\oplus_{\theta}\co^{\theta},\]
where $\theta$ is a central character, 
  $\co^{\theta}\subset\co_{\gtg}$ is a complete
subcategory of $\gtg$-modules 
with generalized central character $\theta$.   One can, therefore,
assume that the functors in question belong to
$Funct(\co^{\theta_{r}},\co^{\theta_{l}})$,
that is, map from $\co^{\theta_{r}}$ to $\co^{\theta_{l}}$.
 Given a weight $\lambda$, 
denote by $M_{\lambda}$ the Verma module 
with highest weight $\lambda$, and $P_{\lambda}$
the indecomposable projective 
module mapping onto the irreducible quotient of $M_{\lambda}$.
Let $\lambda_{r}$ ($\lambda_{l}$ resp.) be a dominant weight such that 
$M_{\lambda_{r}}$ ($M_{\lambda_{l}}$ resp.)
 admits central character $\theta_{r}$ ($\theta_{l}$
resp.) (``Dominant'' here means 
that the corresponding Verma module does not embed in any
Verma module different from itself.) 
From now on it is assumed for simplicity that
 $\lambda_{l},\;\lambda_{r}$ are regular integral. 

\begin{sloppypar}
\begin{theorem}
\label{classprojbg}
There is a bijection
 between the Weyl group $W$ and the isomorphism classes of projective functors
in $Funct(\co^{\theta_{r}},\co^{\theta_{l}})$.
 The bijection is established by assigning
to $w\in W$ the functor
 $\Phi_{w}\in Funct(\co^{\theta_{r}},\co^{\theta_{l}})$, 
such that
$\Phi_{w}(M_{\lambda_{r}})=P_{w\cdot\lambda_{l}}$.
\end{theorem}
\end{sloppypar}
The following result is a key to   Theorem \ref{classprojbg}.

\begin{theorem}

\label{charprojmv}
Any projective functor $\Phi\in Funct(\co^{\theta_{r}},\co^{\theta_{l}})$
 is determined up
to an isomorphism by $\Phi(M_{\lambda_{r}})$.
 The same is true if $\lambda_{r}$ is replaced
with the  antidominant  weight admitting $\theta_{r}$.
\end{theorem}

The way to derive Theorem \ref{classprojbg}.
from Theorem \ref{charprojmv}
is to observe that $\lambda_{r}$ being dominant, $M_{\lambda_{r}}$ is
projective and $V\otimes M_{\lambda_{r}}$ is also. Thus, by
Theorem\ref{charprojmv}, the direct sum
decomposition of $V\otimes ?$ is determined by that of  $V\otimes
M_{\lambda_{r}}$, the latter being delivered by the theory of projectives in
$\co_{\gtg}$.

\begin{remark}
\label{tiltinfindim}
The last sentence of  Theorem \ref{charprojmv} suggests an alternative way to
prove   Theorem \ref{classprojbg}, that is, to consider $V\otimes
M_{\lambda_{r}}$ with an antidominant $\lambda_{r}$. Though not a
projective,  $V\otimes
M_{\lambda_{r}}$   has certain favorable properties summarized by saying
that it is a tilting module. It is this approach which generalizes to the
affine case.
\end{remark}

 Having classified projective functors, it is relatively easy to establish an
equivalence   of (sub)categories of $HCh$ and $\co_{\gtg}$.
Theorem 5.9 of \cite{bernst_gelf} claims that the functor

\begin{equation}
\label{bgisocat1}
HCh(\theta_{l},\theta_{r})\rightarrow \co^{\theta_{l}},\;
H\mapsto \Phi_{H}(M_{\lambda_{r}}),
\end{equation}
is an equivalence of categories. 
Here $HCh(\theta_{l},\theta_{r})$ is a complete subcategory
of  $HCh$ consisting of bimodules 
admitting right central chatacter $\theta_{r}$ and generalized left
central character $\theta_{l}$.

An  important corollary of (\ref{bgisocat1}) is the above mentioned
description of the 2-sided ideal lattice of $U(\gtg)_{\theta}:=
U(\gtg)/U(\gtg)Ker(\theta)$ obtained by Joseph \cite{joseph}.
 Denote by $\Omega(U(\gtg)_{\theta})$ the 
2-sided ideal lattice of $U(\gtg)_{\theta}$ and by $\Omega(M_{\lambda})$
the submodule lattice of $M_{\lambda}$, where $\lambda$ is the dominant
weight related to $\theta$. Then
the map

\begin{equation}
\label{bg_descrofid}
\Omega(U(\gtg)_{\theta})\rightarrow \Omega(M_{\lambda}),
I\mapsto  IM_{\lambda}
\end{equation}
is a lattice equivalence. Indeed,  $U(\gtg)_{\theta}$
is an algebra containing $\gtg$, and hence a $\gtg$-bimodule; its 2-sided ideals as algebra
 are its
submodules as bimodule. Under the equivalence (\ref{bgisocat1}) 
$U(\gtg)_{\theta}$ goes to $M(\lambda)$, because
$U(\gtg)_{\theta}\otimes_{\gtg}M(\lambda)=M(\lambda)$. Thus submodule lattices
of $U(\gtg)_{\theta}$ and $M(\lambda)$ are equivalent. A little extra work
is needed to find the explicit form (\ref{bg_descrofid}) of this 
equivalence.

 The last result  we want to review here belongs to Jantzen
(see \cite{jantz0}, also  \cite{bernst_gelf}) and  establishes
another equivalence of categories based on the notion of  translation
functor. Denote by $\lambda$ the dominant weight lying
in the $W$-orbit of  $\lambda_{1}-
\lambda_{2}$, and by $V_{\lambda}$ the simple $\gtg$-module with highest
weight $\lambda$. (Recall that $\lambda_{l},\;\lambda_{r}$ are supposed
to be integral.)
 For any $\theta$ denote by $p_{\theta}:\co_{\gtg}\rightarrow\co^{\theta}$
the natural projection. Then the functor

\begin{equation}
\label{bf_equiv2}
T_{\theta_{2}}^{\theta_{1}}:\co^{\theta_{2}}\rightarrow \co^{\theta_{1}},\;
T_{\theta_{2}}^{\theta_{1}}(M)=p_{\theta_{1}}(V_{\lambda}\otimes M)
\end{equation}
is an equivalence of categories. The functor $T_{\theta_{2}}^{\theta_{1}}$
is called translation functor.

We finish our review of the semi-simple case by remarking that many results
of \cite{bernst_gelf} are based on, refine and generalize the earlier work,
see e.g. \cite{enr,duflo,vogan,zuck}.

\subsection{An affine analogue.} 
There are many reasons why it is difficult to give
an intelligent definition of a Harish-Chandra bimodule over an affine
Lie algebra $\hgtg$. Some of the difficulties become obvious if one considers universal
enveloping
$U(\hgtg)$ as a model example. Under the diagonal action $U(\hgtg)$    decomposes in a sum
of loop modules. This sum is not direct and multiplicities are infinite. Further, it is
easy to see that the composition series of the
 tensor product of a pair of loop modules always has terms occurring with
 infinite multiplicities. To avoid difficulties of this kind we 
 adopt a functorial point of view.

Thus we are looking for an interesting subcategory in $Funct(\co_{k})$,
$\co_{k}$ being the Bernstein-Gelfand-Gelfand category of $\hgtg$-modules
at level $k$. As an analogue  of the functor $V\otimes ?$
we choose either
\[V_{\lambda}^{k}\tens ?: \co_{k}\rightarrow\co_{k},\; B\mapsto V_{\lambda}^{k}\tens B,\]
where $\tens $ is the Kazhdan-Lusztig tensoring
\cite{kazh_luszt_0,kazh_luszt,kazh_luszt_1,fink},
and $V_{\lambda}^{k}$ is the Weyl module (generalized Verma module in another
terminology) induced in a standard way from the finite dimensional $\gtg$-module
$V_{\lambda}$, or, more generally,
\[A\tens ?: \co_{k}\rightarrow\co_{k},\; B\mapsto A\tens B,\]
where $A$ has a filtration by Weyl modules. To be more precise, we shall consider the following
two cases: (i) $k+h^{\vee} > 0$, and (ii) $k+h^{\vee} < 0$, $h^{\vee} $ being the dual Coxeter
number.

\bigskip

{\bf The case when $k+h^{\vee} > 0$.}
 In this case we confine ourselves to the full
 subcategory $\tco\subset \co_{k}$
 consisting of modules semi-simple with respect to
$\gtg\subset\hgtg$.

The Kazhdan-Lusztig tensoring is
 a subtle thing and many obvious properties of 
$V\otimes ?$ are hard to carry over to the case of $V_{\lambda}^{k}\tens ?$.
For example, the functor $V_{\lambda}^{k}\tens ?$ does not seem to be
exact in general.
There is, however, a case when the analogy is precise -- the affine version
of  translation functor. By \cite{deodgabbkac,rcw},
 there is a direct sum decomposition
\[\tco=\oplus_{(\lambda,k)\in P^{+}_{k}}\tco^{\lambda},\]
and thus a projection
\[p_{\lambda}:\tco\rightarrow\tco^{\lambda},\]
where $P^{+}_{k}$ is the set of
 dominant weights at level $k+h^{\vee}\in\nq_{>}$. (This is
an analogue of the central 
character decomposition for $\gtg$.) We can therefore
define an affine translation functor
\[T_{\mu}^{\lambda}:\tco^{\mu}\rightarrow \tco^{\lambda},\]
by adjusting definition (\ref{bf_equiv2})
 to the affine case (most notably by replacing
$\otimes$ with $\tens$ and
 the finite dimensional $\gtg$-module with an appropriate
Weyl module, for details see \ref {defoftrfunctor}). This construction
was first proposed in \cite{fink} in the case of  negative level 
($k+h^{\vee}< 0$)
representations, but its meaningfulness in our situation is not quite
obvious.

The basic properties of affine translation functors are collected
in Proposition \ref{conseqreltochar}. They are summarized by 
saying that a Weyl module with a dominant highest weight is rigid
and the functor of Kazhdan-Lusztig tensoring with such a module
is exact.  These properties easily imply that
$T_{\mu}^{\lambda}:\tco^{\mu}\rightarrow\tco^{\lambda}$ is
an equivalence of categories (cf. (\ref{bf_equiv2})). This 
theorem refines
results of \cite{deodgabbkac}, where a different version
of translation functors was
defined (in the framework of a general symmetrizable Kac-Moody
algebra) by using the standard tensoring with an integrable module.

The study of Kazhdan-Lusztig tensoring is not easy but rewarding.
A simple translation of Proposition \ref{conseqreltochar} in
the language of vertex operator algebras ( VOA),
see \ref{Vertexoperatorsandvertexoperatoralgebras},
\ref{and_vertex_operator_algebras},
gives the following
affine analogue of the equivalence (\ref{bg_descrofid}). Recall
that by \cite{frzhu} there is a VOA,
$(V_{0}^{k},Y(.,t))$, attached to  $\hgtg$. The Fourier components
of the fields $Y(v,t),v\in V_{0}^{k}$ span a Lie algebra,
$U(\hgtg)_{loc}$. We prove (Theorem \ref{theoronannideals}) that
the ideal lattice of $U(\hgtg)_{loc}$ in the sense of VOA is equivalent to
the submodule lattice of the Weyl module $V_{\lambda}^{k}$
with a dominant highest weight $(\lambda,k),\;k+h^{\vee}\in\nq_{>}$.
Observe that the crucial difference between this statement
and (\ref{bg_descrofid}) is that the associative algebra
$U(\gtg)_{\theta}$ is replaced with a huge Lie algebra
 $U(\hgtg)_{loc}$. Theorem \ref{theoronannideals} generalizes
and refines the well-known result that Fourier components of
the field $e_{\theta}(t)^{k+1}$ annihilate all integrable
modules at a positive level $k$; here $e_{\theta}\in\gtg$ is
a highest root vector.

\bigskip

{\bf The case when $k+h^{\vee} < 0$.} In this case we propose the following affine analogue
of the notion of projective functor. When $k+h^{\vee}<0$, the   theory of projectives in
$\tco$ and $ \co_{k}$ becomes less convenient than its positive level counterpart. Instead of
projectives, we shall use {\em tilting modules}. These can be defined for both
$\tco$ and $\co_{k}$ and their formal characters linearly span the Grothendieck rings of
each of the categories. To give a couple of examples, remark that
a Verma module with
an antidominant highest weight is  a tilting object of
$\co_{k}$; likewise, a Weyl module with
an antidominant highest weight is  a tilting object of
$\tco$. (It is worth mentioning that neither $\co_{k}$ nor $\tco$ has
 modules with
dominant highest weight.) 

Each tilting module is a direct sum of indecomposable ones and the isomorphism
classes of the latter  are in 1-1 correspondence with highest weights.

We define a {\em tilting functor} to be a direct summand of the functor $p_{\lambda_{l}}
\circ(W\tens ?)\in Funct(^{fam}\co_{k}^{\lambda_{r}},^{fam}\co_{k}^{\lambda_{l}})$, where
$W\in \tco$ is tilting. Here $^{fam}\co^{\lambda}_{k}$ is a complete
subcategory of $\co^{\lambda}_{k}$ consisting of modules which can be
included in a family analytically depending on $k$. We prove (see Theorem \ref{class_titl_fun})
that isomorphism classes of indecomposable tilting
functors are in 1-1 correspondence with elements of affine Weyl group $W_{k}$ under
the asumption that all the highest weights in question are regular. Further,
as we mentioned in Remark \ref{tiltinfindim},
the correspondence is established in the same  way as the 
one used in Theorem \ref{classprojbg} except that dominant highest weights are  replaced
with antidominant ones and projective modules are  replaced with tilting
ones.
	It follows that the category of tilting functors is
 semisimple and is isomorphic to the
subcategory of tilting modules in $\co_{k}^{\lambda_{l}}$. We remark that 
 Theorem \ref{class_titl_fun} is derived from Theorem \ref{titlfunvsitsvalue}
in the way analogous to the
one used to derive Theorem \ref{classprojbg} from Theorem \ref{charprojmv}. However, when
proving Theorem \ref{titlfunvsitsvalue}, we do not have have at our disposal
   the important tool used by Bernstein and S.Gelfand, namely, the bimodule interpretation of
projective functors, in particular the universal enveloping algebra and Kostant's theorem.
Instead, we rely on some techniques borrowed from conformal field theory.

Many of the difficulties which one encounters when working with $\tens$ were removed
in \cite{kazh_luszt,kazh_luszt_1}. For instance, the functor $W\tens?$ is known to be
exact and, therefore, induces a homomorphism of the Grothendieck rings. We show that this
homomorphism is  multiplication by a certain element of the group algebra of $W_{k}$.
In the case when $\lambda_{r}=\lambda_{l}$, the category of tilting functors is
naturally a monoidal category, the tensor product being simply the 
composition of functors,
and the previous sentence implies that the Grothendieck ring of this category is isomorphic to
the group algebra of $W_{k}$, see Theorem \ref{grothringoftilt}.

 \begin{sloppypar}
{\bf Acknowledgements.} We thank P.Etingof, B.Feigin and
G.Zuckerman for
illuminating discussions and to I.Mircovi\v{c} who told us about
his unpublished work \cite{mirc} containing results reminiscent of some of
ours.
Part of the work was done when we were visiting E.Shr\"{o}dinger Institute
in Vienna in June 1996 . We are grateful to the organizers of the conference
on Representation  Theory and Applications to Mathematical Physics for
invitation and to the Institute for support. The second author enjoyed
the support of IHES and discussions with K.Gawedzki during July-August of
1996. In the present form this paper was first reported in Feigin's
seminar at the Moscow Independent University.
\end{sloppypar}

\section{preliminaries}
\label{preliminaries}
\subsection{ }
\label{prelim_reps}
  The following is a list of essentials which
will be used but will not be explained.

$\gtg$ is a simple finite dimensional Lie algebra
with a fixed triangular decomposition;
in particular with a fixed Cartan subalgebra
$\gth\subset\gtg$; unless otherwise mentioned the dual
$\gth^{\ast}$ is understood as the real part
 $\gth^{\ast}_{\bf R}$;

the action ($\lambda\mapsto 
w\lambda$) and the shifted (by $\rho$)
 action ($\lambda\mapsto w\cdot\lambda$) 
of the Weyl group $W$ on $\gth^{\ast}$ preserving the
weight lattice $P\in\gth^{\ast}$; denote by $\bar{C}$ the Weyl chamber --
a fundamental domain for the shifted action attached to the fixed
triangular decomposition; $P^{+}=P\cap C$ where $C\subset\bar{C}$
is the interior;

the $\co_{\gtg}$ category of $\gtg$-modules attached to
the triangular decomposition;

a Verma module $M_{\lambda}\in\co_{\gtg},\; \lambda
\in\gth$ and a simple finite dimensional 
module
 $V_{\lambda},\; \lambda\in P^{+}\subset P$;

the affine Lie algebra $\hgtg=\gtg\otimes\nc((z))\oplus\nc K$
 and the ``generalized''
Borel subalgebra $\hgtg_{\geq}=\gtg\otimes\nc[[z]]\oplus\nc K$;

$\co_{k}$ -- the category of $\hgtg$-modules at level $k$
 (i.e. $K\mapsto k$),
and the full subcategory $\tco\subset\co_{k}$ consisting of 
$\hgtg$-modules semisimple over $\gtg\subset\hgtg$;

$M_{\lambda}^{k}=Ind_{\hgtg_{\geq}}^{\hgtg}M_{\lambda}
\in\co_{k},\; \lambda\in\gth^{\ast}$ is a Verma module;
 $V_{\lambda}^{k}=Ind_{\hgtg_{\geq}}^{\hgtg}V_{\lambda}
\in\tco,\;\lambda\in P^{+}$ is a Weyl module; more generally,
if $V\in\co_{\gtg}$ is a $\gtg$-module, then $V^{k}\in\co_{k}$
 is a $\hgtg$-module
obtained by inducing from $V$; obviously, $V_{\lambda}^{k}$ is a
quotient of $M_{\lambda}^{k}$; each simple module in $\co_{k}$
 is a quotient
of $M_{\lambda}^{k}$ for some $\lambda,\; k$; denote this
module by $L_{\lambda}^{k}$;

if $k\not{\in}\nq$, then $\tco$ is semi-simple, each object being
a direct sum of Weyl modules; the analogue of this
statement for $\co_{k}$ is the equivalence $\co_{k}\approx \co_{\gtg}$
obtained by Finkelberg \cite{fink};

for $k+h^{\vee}=p/q\in\nq_{>}$ consider the affine Weyl group
$W_{k}=pQ\propto W$, where $Q$ is a root lattice of $\gtg$; there
is the usual and the dotted (shifted) action of $W_{k}$ on
$\gth^{\ast}$; the fundamental domain for the latter is
$\bar{C}_{aff}=\bar{C}\cap\{\lambda: 0\leq (\lambda+\rho,\theta)\leq p\}$,
 where
$\theta$ is the highest root of $\gtg$; set
$P^{+}_{k}=P^{+}\cap C_{aff}$ where $C_{aff}\subset \bar{C}_{aff}$ is the
interior; at some points it is important to ensure that $P^{+}_{k}$ contains
at least one non-zero weight; for that, in the case of some exceptional
root systems, $k$ should be sufficiently large (see, for example,
\cite{fink} sect. 2.6), and this will be our assumption throughout
the text;
 call $\lambda\in P^{+}_{k}$ (sometimes
$(\lambda,k)$  if $\lambda$ satisfies this condition) dominant;
if  $k+h^{\vee}=p/q\in\nq_{<}$, one defines $W_{k}$ and an
antidominant weight in a similar way;

by \cite{deodgabbkac,rcw}, $\co_{k}=\oplus_{\lambda\in P^{+}_{k}}
\co_{k}^{\lambda}$, where $\co_{k}^{\lambda}$ is a full subcategory
consisting of modules whose composition series contain only
irreducible modules $L_{w\cdot\lambda}^{k}, w\in W_{k}$; a similar
decomposition is true for $\tco$.

\bigskip

{\em Duality Functors.} Given a vector space $W$, denote by $W^{d}$ its
total dual. If $W$ is a Lie algebra module, then so is $W^{d}$.

 Given a vector space $W$ carrying  a gradation by finite
dimensional subspaces, denote by $D(W)$ its restricted dual.

Objects of $\tco$ are canonically graded. Denote by $D:\tco\rightarrow
\tco$, $M\mapsto D(M)$ the functor such that the $\hgtg$-module structure
is defined by precomposing the canonical action on the dual space with
an automorphism $\hgtg\rightarrow\hgtg$, 
$g\otimes z^{n}\mapsto g\otimes(-z)^{-n}$. In a similar manner one defines the duality
$D:\co_{k}\rightarrow\co_{k}$.

The functors $^{d}, D(.)$ are exact.

There is an involution $\bar{ }: P^{+}\rightarrow P^{+}$ so that
$V_{\lambda}^{d}=V_{\bar{\lambda}}$.

\subsection{Two lemmas on geometry of weights}

\subsubsection{ }
The following is proved in \cite{jantz} Lemma 7.7.
\begin{lemma}
\label{lemmaongeomofw}
Suppose:

(i) $(\lambda,k),(\mu,k)\in P^{+}_{k}$ are regular;

(ii)$\bar{w}\in W$ satifies $\bar{w}(\lambda-\mu)\in P^{+}$;

(iii) $\nu$ is a weight of $V_{\bar{w}(\lambda-\mu)}$ such that
$w_{1}\cdot\lambda=w\cdot\mu +\nu$ for some $w,w_{1}\in W_{k}$.

Then: $w_{1}=w$ and $\nu\in W(\lambda-\mu)$.
\end{lemma}

\subsubsection{ }
\label{defofordandlem}
The Bruhat ordering on $W_{k}$ determines a partial ordering on
$\gth^{\ast}$: $\mu<_{k}\nu$ if and only if there is an (anti)dominant
$\lambda$ so that $\mu=w_{1}\cdot\lambda,\;\nu=w_{2}\cdot\lambda$ for
some $w_{1}<w_{2}\in W_{k}$.

\begin{lemma}
\label{longowbg}
Let $\lambda$ be a regular antidominant weight. If
$\lambda+\psi<_{k}\lambda+\phi$, then $|\phi|>|\psi|$, where $|.|$
is the length function
coming from the canonical inner product on $\gth^{\ast}$.
\end{lemma}

{\bf Proof} repeats word for word the proof of Lemma 1.5 in Appendix 1 of
\cite{bernst_gelf} except that instead of hyperplanes one has to consider
affine hyperplanes.

\section{ the kazhdan-lusztig tensoring}
\label{kazhdluszttens}

Kazhdan and Lusztig \cite{kazh_luszt_0,kazh_luszt,kazh_luszt_1} 
(inspired by Drinfeld \cite{drinf_1}) defined
a covariant bifunctor

\begin{equation}
\label{notattens}
\tco\times\tco\rightarrow\tco,\; A,B\mapsto A\tens B.
\end{equation}

We shall review its definition and main properties.
\subsection {Definition}
\label{twodefinitions}

\subsubsection{ The set-up}
\label{thesetup}

The notation to be used is as follows:

 $z$ is a  once and for all fixed coordinate on $\cp$;

$L\gtg^{P}, P\in\cp$ is the loop algebra attached to $P$; in other words,
$L\gtg^{P}=\gtg\otimes\nc((z-P)), P\in\nc$, and 
$L\gtg^{\infty}=\gtg\otimes\nc((z^{-1}))$;

more generally, if $P=\{P_{1},...,P_{m}\}\subset\cp$, then 
 \[L\gtg^{P}=\oplus_{i=1}^{m} L\gtg^{P_{i}};\] 

$\hgtg^{P}=L\gtg^{P}\oplus \nc K, P\in\cp$ is the affine algebra attached to the point $P$ --
the canonical central extension of $L\gtg^{P}$; of course, $\hgtg^0=\hgtg$;

more generally, if $P=\{P_{1},...,P_{m}\}\subset\cp$, then 
 $\hgtg^{P}$ is the direct sum of  $\hgtg^{P_{i}},\; i=1,...,m$ modulo the
relation: all canonical central elements $K$ (one in each copy) are equal
each other;

$\Gamma=\gtg\otimes\nc[z,z^{-1},(z-1)^{-1}]$; $\Gamma$ is obviously a Lie
algebra.

The Laurent series expansions at points $\infty,1,0$ produce the Lie algebra
homomorphism
\[\epsilon: \Gamma\rightarrow L\gtg^{\{\infty,1,0\}}.\]

\begin{lemma}
\label{trehtochech}
The map $\epsilon$ lifts to a Lie algebra homomorphism
\[ \Gamma\rightarrow \hgtg^{\{\infty,1,0\}}.\]
\end{lemma}

Proof consists of using the residue theorem, see \cite{kazh_luszt}.

By pull-back, any $\hgtg^{\{\infty,1,0\}}$-module is canonically a
$\Gamma$-module.
 Further, any $A\in\tco$ is canonically a $\hgtg^{P}$-module
for any $P$  -- by the obvious change of coordinates; refer to this
as attaching $A$ to $P\in\cp$. Given $A,B,C\in\tco$, we shall regard $A\otimes B\otimes C$
as a  $\hgtg^{\{\infty,1,0\}}$-module meaning that  $\hgtg^{\infty}$ acts on
$A$,  $\hgtg^{1}$ on $B$,  $\hgtg^{0}$ on $C$. (There is an obvious
ambiguity in this notation.) There arises the space of coinvariants
 \[(A\otimes B\otimes C)_{\Gamma}= (A\otimes B\otimes C)/\Gamma(A\otimes B\otimes C).\]

This construction easily generalizes to the case when instead of three points
-- $\infty,1,0$ -- there are $m$ points,   $m$ modules and instead of $\Gamma$
one considers the Lie algebra of rational functions on $\cp$ with $m$ punctures
with values in $\gtg$. We shall be mostly interested in the case $m=3$
and sometimes
in the case $m=2$. If $m=2$, then $\Gamma$ becomes
 $\tgtg=\gtg\otimes\nc[z,z^{-1}]$.

\begin{lemma}
\label{morphbetw2andcoinv}
Suppose $D(B)$ is attached to $\infty$, $A$ to 0. Then
\[Hom_{\hgtg}(A,B)=((D(B)\otimes A)_{\tgtg})^{d}.\]
\end{lemma}

Proof can be found in \cite{kazh_luszt}; the reader may also observe that
the arguments from \ref{morphismsandcoinvariants} are easily adjusted to this
case.

\subsubsection{Definition}
\label{Theseconddefinition}
Let $\hg$ be the central extension of $\Gamma$, the cocycle being defined
as usual except that one takes the sum of residues at $\infty$ and 1. Let
$\Gamma(0)\subset\hg$ be the subalgebra consisting of functions vanishing
at 0. Obviously, $\Gamma(0)$ can also be regarded as a subalgebra of $\Gamma$.

Consider the (total) dual space  $(A\otimes B)^{d}$; it is naturally a
 $\hg$-module.  $(A\otimes B)^{d}$ carries the increasing filtration
$\{ (A\otimes B)^{d}(N)\}$, where

\begin{multline}
(A\otimes B)^{d}(N)\\
=\{x\in (A\otimes B)^{d}: \gamma_{1}\cdots\gamma_{N}x=0 \mbox{ if all } \gamma_{i}\in \Gamma(0),
x\in (A\otimes B)\}.
\end{multline}

The space 
$\cup_{N\geq 1} (A\otimes B)^{d}(N)$ is naturally a $\hgtg$-module. The passage
from $(A\otimes B)^{d}$ to $\cup_{N\geq 1} (A\otimes B)^{d}(N)$ (or its obvious
versions) is  often called  a functor of smooth vectors.

{\bf Define} 
\begin{equation}
\label{seconddefform}
A\tens B=D(\bigcup_{N\geq 1} (A\otimes B)^{d}(N)).
\end{equation}

\begin{lemma}
\label{rightexactness}
The functor $\tens : \tco\times\tco\rightarrow \tco$ is right exact in
each variable.
\end{lemma}

{\bf Proof}  (see {\em loc. cit.}) The functor $\tens$ is a composition of two dualizations, $^d$ and $D(.)$,
and the functor of smooth vectors. It is enough to remark that the first two
are exact while the last is only left exact.$\qed$

\subsection {Some properties of $\tens$.}

\subsubsection{ }
For the future reference we collect some of the properties
of $\tens$ in the following

\begin{theorem}
\label{mainthonklus}

(i)
\[Hom_{\hgtg}(A\tens B, D(C))= Hom_{\hgtg}(C,D(A\tens B))
=((A\otimes B\otimes C)_{\Gamma})^{d}.\]

(ii) If $A,B\in\tco$ have a Weyl filtration, then $A\tens B$ has
also. (Here by Weyl filtration we mean a filtration such that its quotients
are Weyl modules.)

(iii) If $k\not\in\nq$, then $V_{\lambda}^{k}\tens V_{\mu}^{k}
=(V_{\lambda}\otimes V_{\mu})^{k}$.

(iv) For any $k\in\nc$, $V_{\lambda}^{k}\tens V_{\mu}^{k}$ has
 a Weyl filtration (see(ii)), the multiplicity of $V_{\nu}^{k}$
being equal $(V_{\lambda}\otimes V_{\mu}: V_{\nu})$ (c.f. (iii)).

(v) There is an isomorphism $A\tens V_{0}^{k}\rightarrow A$
for any $A\in\tco$.

(vi) There are commutativity and associativity morphisms
$A\tens B\approx B\tens A$ and $(A\tens B)\tens C\approx A\tens (B\tens C)$
which endow $\tco$ with the structure of a braided monoidal category.

\end{theorem}

\subsubsection{Morphisms and coinvariants.}
\label{morphismsandcoinvariants}
 The description
of morphisms in terms of coinvariants (see Theorem \ref{mainthonklus}(i))
is the hallmark of this theory. Let us briefly explain why (i) holds.
There is the obvious isomorphism of vector spaces
\[(A\otimes B\otimes C)^{d}\rightarrow  Hom_{\nc}(C,(A\otimes B)^{d}).\]
It induces the map

\[((A\otimes B\otimes C)_{\Gamma}^{d}\rightarrow  Hom_{\hg}(C,(A\otimes B)^{d}).\]
By $\hg$-linearity, it actually gives the map

\[((A\otimes B\otimes C)_{\Gamma}^{d}\rightarrow  Hom_{\hg}(C,\bigcup_{N\geq 1} 
(A\otimes B)^{d}(N)).\]
It remains to look at (\ref{seconddefform}) and note that $\hg$ is dense in $\hgtg$.

\subsubsection{Using the spaces of coinvariants.}
\label{Usingthespacesofcoinvariants}
 A lot about the functor $\tens$ easily
follows from Theorem \ref{mainthonklus}(i). As an example, let us derive  (v). By (i),

\[Hom_{\hgtg}(A\tens V_{0}^{k},B)=((A\otimes V_{0}^{k}\otimes D(B))_{\Gamma})^{d}
\mbox{ for any } B\in\tco.\]

As $V_{0}^{k}=Ind_{\hgtg_{\geq}}^{\hgtg}\nc$, the Frobenius reciprocity
gives

\[(A\otimes V_{0}^{k}\otimes B)_{\Gamma}=(A\otimes D(B))_{\tgtg},\]
the latter space being  $Hom_{\hgtg}(A,B)$ by Lemma \ref{morphbetw2andcoinv}. 
We see that the spaces
of morphisms of the modules $A$ and $A\tens V_{0}^{k}$ are equal, hence
so are the modules.

Replacing in this argument $\nc$ with a suitable finite dimensional
$\gtg$-module and repeating it three times one gets

\begin{equation}
\label{morphofweylintodual}
Hom_{\hgtg}(V_{\lambda}^{k}\tens V_{\mu}^{k}, D(V_{\nu}^{k}))
=Hom_{\gtg}(V_{\lambda}\otimes V_{\mu}, V_{\bar{\nu}}).
\end{equation}

As for generic $k$ $D(V_{\nu}^{k})\approx V_{\nu}^{k}$
(see \ref{prelim_reps}), (\ref{morphofweylintodual}) along with 
Theorem \ref{mainthonklus}(i) implies Theorem \ref{mainthonklus}(iii).

\section{Affine translation functors}

\subsection{Definition.}
\label{defoftrfunctor}
For any  $(\lambda,k )\in P^{+}_{k}$ denote by $\tco^{\lambda}$ the full
subcategory  
of $\tco$ consisting of modules whose composition factors all have highest weights
 lying in the orbit $W_{k}\cdot(\lambda,k)$. There arises the projection
\[p_{\lambda}:\tco\rightarrow \tco^{\lambda}.\]
This all has been reviewed in \ref{prelim_reps}.

Given  $(\lambda,k ),(\mu,k)\in P^{+}_{k}$, pick $\bar{w}\in W$ so that
$\bar{w}(\lambda-\mu)\in P^{+}$. It is easy to see that then
$(\bar{w}(\lambda-\mu),k)\in P^{+}_{k}$.

{\bf Define } the translation functor

\begin{equation}
\label{defoftrfunctor_form}
T_{\mu}^{\lambda}:\; \tco^{\mu}\rightarrow \tco^{\lambda}
A\mapsto p_{\lambda}(V_{\bar{w}(\lambda-\mu)}^{k}\tens A).
\end{equation}
This functor was first introduced by Finkelberg \cite{fink} who, however, considered it only
for $k<0$. 

As an immediate corollary of the definition, one has
\begin{equation}
\label{howtoadjointtlm}
T_{\lambda}^{\mu}=p_{\mu}\circ ((V_{\bar{w}(\lambda-\mu)}^{d})^{k}\tens ?)
\end{equation}

\subsection{Rigidity of  Weyl modules with dominant highest weight.}
\label{rigidityofweylmod}
\begin{lemma}
\label{actoftlmonweyls}

 If $(\lambda,k),(\mu,k)$ are regular (i.e. off the affine walls) and $w\in W_{k}$
satisfies $w\cdot\mu\in P^{+}$, then
\[T_{\mu}^{\lambda}(V_{w\cdot\mu}^{k})=V_{w\cdot\lambda}^{k}.\]

\end{lemma}

{\bf Proof.}  By Theorem \ref{mainthonklus} (iv), $T_{\mu}^{\lambda}(V_{w\cdot\mu}^{k})$
has a filtration with quotients isomorphic to $V_{w_{1}\cdot\lambda}^{k}$,
$w_{1}\in W_{k}$ such that $w_{1}\cdot\lambda =w\cdot\mu+\nu$, $\nu$ being
a weight of $V_{\bar{w}(\lambda-\mu)}$. By Lemma \ref{lemmaongeomofw}, $w_{1}=w$. This implies that this filtration has only one term, 
$V_{w\cdot\lambda}^{k}$. $\qed$

\begin{corollary}
\label{rigofvlk}
If $(\lambda,k)\in P^{+}_{k}$ is regular, then $V_{0}^{k}$ is a direct summand of
$V_{\lambda}^{k}\tens V_{\bar{\lambda}}^{k}$.
\end{corollary}

{\bf Proof.} Of course $(0,k)$ is dominant regular
and $p_{0}A$ is a direct summand of $A$. It remains
to observe that $T_{\bar{\lambda}}^{0}V_{\bar{\lambda}}^{k}=
p_{0}(V_{\lambda}^{k}\tens V_{\bar{\lambda}}^{k})$ and use 
 Lemma \ref{actoftlmonweyls} to get
 $T_{\bar{\lambda}}^{0}V_{\bar{\lambda}}^{k}=
V_{0}^{k}$ . $\qed$

We get the maps
\[i_{\lambda}: V_{0}^{k}\rightarrow V_{\lambda}^{k}\tens
 V_{\bar{\lambda}}^{k},\;
e_{\lambda}: V_{\bar{\lambda}}^{k}\tens V_{\lambda}^{k}\rightarrow V_{0}^{k}.\]

Observing that the maps between $\tens$-products of Weyl
modules are uniquely determined by the induced maps of
the corresponding finite dimensional $\gtg$-modules (Theorem \ref{mainthonklus} and 
(\ref{morphofweylintodual}) ), we see that
 we can normalize $i_{\lambda},e_{\lambda}$ so that the compositions

\begin{multline}
V_{\lambda}^{k}=V_{0}^{k}\tens V_{\lambda}^{k}\stackrel{i_{\lambda}\otimes id}{\rightarrow}
V_{\lambda}^{k}\tens V_{\bar{\lambda}}^{k}\tens V_{\lambda}^{k}
\stackrel{id\otimes e_{\lambda}}{\rightarrow} V_{\lambda}^{k}\\
V_{\bar{\lambda}}^{k}=V_{\bar{\lambda}}^{k}\tens V_{0}^{k}\stackrel
{id\otimes i_{\lambda}}
{\rightarrow}
V_{\bar{\lambda}}^{k}\tens V_{\lambda}^{k}\tens V_{\bar{\lambda}}^{k}
\stackrel{ e_{\lambda}\otimes id}{\rightarrow} V_{\bar{\lambda}}^{k},
\label{feofrigid}
\end{multline}

are equal to the identity. By definition (see e.g. \cite{kazh_luszt_1} III, 
Appendix) we have

\begin{corollary}
\label{rigidity}
If $(\lambda, k)\in P^{+}_{k}$, then $V_{\lambda}^{k}$ and 
$V_{\bar{\lambda}}^{k}$
are rigid.
\end{corollary}

Consider the functor $V_{\lambda}^{k}\tens ?: \tco\rightarrow\tco,
 M\mapsto V_{\lambda}^{k}
\tens M.$

\begin{corollary}
\label{adjointnessproperty}
(i)If $(\lambda, k)\in P^{+}_{k}$, then the functors $V_{\lambda}^{k}\tens ?$
 and $V_{\bar{\lambda}}^{k}\tens ?$ are adjoint, i.e. there is a functor ismorphism

\[Hom_{\hgtg}(V_{\lambda}^{k}\tens A,B)=Hom_{\hgtg}( A,V_{\bar{\lambda}}^{k}\tens B).\]

(ii) If $(\lambda, k)\in P^{+}_{k}$, then the functors $V_{\lambda}^{k}\tens ?$
 and $V_{\bar{\lambda}}^{k}\tens ?$ are exact, i.e. send exact short sequences to
exact ones.
\end{corollary}

{\bf Proof} is standard; for the reader's convenience we reproduce the one from  \cite{kazh_luszt_1} III, 
Appendix.
To prove (i), consider two composition maps

\[\phi:Hom_{\hgtg}(V_{\lambda}^{k}\tens A,B)\rightarrow 
Hom_{\hgtg}(V_{\bar{\lambda}}^{k}\tens  V_{\lambda}^{k}\tens A,V_{\bar{\lambda}}^{k}\tens B)
\stackrel{i_{\bar{\lambda}}}{\rightarrow} 
Hom_{\hgtg}( A,V_{\bar{\lambda}}^{k}\tens B),\]

\[\psi:Hom_{\hgtg}( A,V_{\bar{\lambda}}^{k}\tens B)\rightarrow 
Hom_{\hgtg}(V_{\lambda}^{k}\tens A,V_{\lambda}^{k}\tens  V_{\bar{\lambda}}^{k}\tens B)
\stackrel{e_{\bar{\lambda}}}{\rightarrow }
Hom_{\hgtg}(V_{\lambda}^{k}\tens A, B).\]

By (\ref{feofrigid}), the compositions $\phi\circ\psi$ and $\psi\circ\phi$ are equal to
the identity.

(ii) is an easy consequence of (i): we have to prove that $B_{1}\rightarrow B_{2}$ is
a monomorphism implies that $V_{\lambda}^{k}\tens B_{1}\rightarrow 
V_{\lambda}^{k}\tens B_{2}$ is also, or, equivalently, that for any $A\in\tco$
the induced map
\[Hom_{\hgtg}(A,V_{\lambda}^{k}\tens B_{1})\rightarrow 
Hom_{\hgtg}(A,V_{\lambda}^{k}\tens  B_{2})\]
is also a monomorphism. By (i), it is equivalent to proving that 
\[Hom_{\hgtg}(V_{\bar{\lambda}}^{k}\tens A, B_{1})\rightarrow 
Hom_{\hgtg}(V_{\bar{\lambda}}^{k}\tens A, B_{2})\]
is a monomorphism, but this is an obvious corollary of injectivity of the map
$B_{1}\rightarrow B_{2}$. $\qed$

\subsection{Properties of affine translation functors.}
Recall that there is the notion of  formal character $ch A$ for
any $A\in\tco^{\lambda}$, see e.g. \cite{deodgabbkac}. There arises an abelian
group of characters, each of the following sets being a topological basis in it:

$\{ch V_{w\cdot\lambda}^{k},\; w\in W_{k}\}$,
$\{ch L_{w\cdot\lambda}^{k},\; w\in W_{k}\}$.
Of course the symbols $ch V_{w\cdot\lambda}^{k},\; ch L_{w\cdot\lambda}^{k}$
should be ignored unless $w\cdot\lambda\in P^{+}$. Observe that 
\begin{equation}
\label{mezhduvil}
ch A= \sum_{w\geq w_{o}}\bar{n}_{w}ch L_{w\cdot\mu}^{k}\Leftrightarrow
ch A= \sum_{w\geq w_{o}}n_{w}ch V_{w\cdot\mu}^{k}
\end{equation}

\begin{proposition}
\label{conseqreltochar}
Let  $(\lambda,k),(\mu,k)$ be regular dominant.

(i) $T_{\mu}^{\lambda}$ is exact;

(ii)  $T_{\mu}^{\lambda}, T_{\lambda}^{\mu}$ are adjoint to each other;

(iii)  If $ch A=\sum_{w\in W_{k}}n_{w}ch V_{w\cdot\mu}^{k}$, then
 $ch T_{\mu}^{\lambda} A=\sum_{w\in W_{k}}n_{w}ch V_{w\cdot\lambda}^{k}$.

(iv) $T_{\mu}^{\lambda}(L_{w\cdot\mu}^{k})=L_{w\cdot\lambda}^{k}$;

(v) More generally,  $T_{\mu}^{\lambda}(.)$ establishes
an equivalence of the submodule lattices of 
$V_{w\cdot\mu}^{k}$ and $V_{w\cdot\lambda}^{k}$.
\end{proposition} 

{\bf Proof.}
(i)  $T_{\mu}^{\lambda}$ is exact as a composition of the exact functors
$p_{\lambda}$ and $(V_{\bar{w}(\lambda-\mu)}^{d})^{k}\tens ?$, see
Corollary \ref{adjointnessproperty} (ii).

(ii) By Corollary \ref{adjointnessproperty} (i),
one has for any $A\in\tco^{\mu},\; B\in\tco^{\lambda}$

\begin{multline}
Hom_{\hgtg}(T_{\mu}^{\lambda}A,B)=
Hom_{\hgtg}(p_{\lambda}(V_{\bar{w}(\lambda-\mu)}^{k}\tens A),B)
=Hom_{\hgtg}(V_{\bar{w}(\lambda-\mu)}^{k}\tens A,B)\\=
Hom_{\hgtg}( A,(V_{\bar{w}(\lambda-\mu)}^{d})^{k}\tens B)
=Hom_{\hgtg}( A,p_{\mu}(V_{\bar{w}(\lambda-\mu)}^{d})^{k}\tens B)\\
=Hom_{\hgtg}( A, T_{\lambda}^{\mu}B).
\end{multline}

(iii) follows at once from (i) ( if one uses the local composition series, see e.g. \cite{deodgabbkac}).

(iv) Let $T_{\mu}^{\lambda}(L_{w_{0}\cdot\mu}^{k})$ be reducible. There
arises an exact sequence with non-zero $N$
\[0\rightarrow N\rightarrow  T_{\mu}^{\lambda}(L_{w_{0}\cdot\mu}^{k})
\rightarrow L_{w_{0}\cdot\lambda}^{k}\rightarrow 0.\]

Applying $T_{\lambda}^{\mu}$ to it one gets

\[0\rightarrow T_{\lambda}^{\mu}( N)\rightarrow  T_{\lambda}^{\mu}(T_{\mu}^{\lambda}(L_{w_{0}\cdot\mu}^{k}))
\rightarrow T_{\lambda}^{\mu}(L_{w_{0}\cdot\lambda}^{k})\rightarrow 0.\]

By (iii) and (\ref{mezhduvil}),
 $ch(T_{\lambda}^{\mu}(T_{\mu}^{\lambda}(L_{w_{0}\cdot\mu}^{k})))=
ch L_{w_{0}\cdot\mu}^{k}$ and $ch T_{\lambda}^{\mu}( N)\neq 0$; therefore 
$ch T_{\lambda}^{\mu}(L_{w_{0}\cdot\lambda}^{k})< ch L_{w_{0}\cdot\mu}^{k}$.
Contradiction. 

(v) Here proof is an obvious version of that of (iv). By
  (ii) it is enough to show that 
if $A\subset B\subset V_{w\cdot\mu}^{k}$, then 
 $T_{\mu}^{\lambda}(A)\subset T_{\mu}^{\lambda}( B)\subset V_{w\lambda}^{k}$.
Using (\ref{mezhduvil}) and passing to quotients, if necessary, the
problem is reduced to the case when $B$ is a highest weight module. In this
case the arguments of (ii) go through practically unchanged.
$\qed$

\subsection { }

\begin{theorem}
\label{equivofcat}
The functor $T_{\mu}^{\lambda}:\tco^{\mu}\rightarrow\tco^{\lambda}$ is
an equivalence of categories.
\end{theorem}

{\bf Proof.} It is enough show that $T_{\mu}^{\lambda}\circ T_{\lambda}^{\mu}:
\tco^{\lambda}\rightarrow \tco^{\lambda}$ is equivalent to the identity.
In other words, we want to show that $id: A\rightarrow A, \; A\in\tco^{\lambda}$ is 
transformed into an isomorphism in
 $Hom_{\hgtg}(T_{\mu}^{\lambda}\circ T_{\lambda}^{\mu}(A), A)$. We already know this when $A$
 is simple, see
Corollary \ref{conseqreltochar} (ii). Using the local composition
series one proves it for an arbitrary $A$.

An alternative way to prove the theorem is to observe that by Corollary
\ref{rigofvlk} the action of  $T_{\mu}^{\lambda}\circ T_{\lambda}^{\mu}$ is
equivalent to that of $V_{0}^{k}\tens ?$, the latter being equivalent
to $id$ by Theorem \ref{mainthonklus} (v). $\qed$

\subsection { Generalizing from $\tco$ to $\co_{k}$.}
\label{genfromtcotook}

Our two key results -- Proposition \ref{conseqreltochar} and Theorem \ref
{equivofcat} -- can be carried over to the category $\co_{k}$. Let us briefly
explain it. We will be using subcategories $\co_{k}^{\lambda}\subset
\co_{k}$ (see \ref{prelim_reps}) only when
$k+h^{\vee}\in\nq_{>}$ and $\lambda$ is integral, although the last condition
can be easily relaxed.

It is non-trivial (if at all meaningful ) to carry the Kazhdan-Lusztig tensoring over
 to the entire $\co_{k}$.  It is
however straightforward to extend it to the functor
\[\tens:\;\tco\times\co_{k}\rightarrow \co_{k},\]
as proposed by Finkelberg \cite{fink}.  One basic property of this operation
absolutely analogous (along with the proof) to Theorem \ref{mainthonklus} (iv)
is as follows. 

\begin{lemma}
\label{multinfiltr}

If  $A\in\tco$ has a filtration by Weyl modules and
$B\in\co_{k}$ has a filtration by Verma modules, then
 $A\tens B$ also has a filtration
by Verma modules. Further the multipliciites are the same as
in the finite dimensional case; for instance
\begin{equation}
\label{multinfiltr_eq}
(V_{\lambda}^{k}\tens M_{\mu}^{k}:M_{\nu}^{k})=
(V_{\lambda}\otimes M_{\mu}: M_{\nu}).
\end{equation}
\end{lemma}

Given this one can easily inspect our exposition of affine translation
functors and
observe that quite a lot  carries over to the setting of
$\co_{k}$ word for word except that at the
appropriate places Weyl modules are to be changed for the corresponding
Verma modules. Here are some examples:

(i) definition of
 $T_{\mu}^{\lambda}:\co_{k}^{\mu}\rightarrow\co_{k}^{\lambda}$;

(ii) the Verma filtration of $V_{\lambda}^{k}\tens M_{w\cdot\mu}^{k}, w\in W_{k}$ and
 Lemma \ref{lemmaongeomofw} imply that $T_{\mu}^{\lambda}(M_{w\cdot\mu}^{k})=
M_{w\cdot\lambda}^{k}$ if $(\mu,k),(\lambda,k)$ are regular
 (c.f. Lemma \ref{actoftlmonweyls}); observe that we can now drop the condition that
$w\cdot\mu\in P^{+}$;

(iii) therefore Proposition \ref{conseqreltochar} holds with the indicated changes.

 We get 
\begin{theorem}
\label{eqofcat_forbigg}
The functor $T_{\mu}^{\lambda}:\co^{\mu}_{k}\rightarrow\co^{\lambda}_{k}$ is
an equivalence of categories if 
 $\lambda,\;\mu$ are integral and both belong to the same Weyl chamber.
\end{theorem}

\section {annihilating ideals of highest weight modules}

\subsection{Vertex operators and ...}
 
The usual tensor functor $\otimes: M,N\mapsto M\otimes N$ has
the following fundamental (and trivial) property: there is a natural map
\begin{multline}
\label{fundmapatttoustens}
N\rightarrow Hom_{\nc}(M,M\otimes N)\\
n\mapsto n(.) \mbox{ such that } n(m)=m\otimes n.
\end{multline}

Here  we shall explain the $\tens$-analogue of this map

\subsubsection{ }
\label{Vertexoperatorsandvertexoperatoralgebras}

By Theorem \ref{mainthonklus} (v), $A\tens V_{0}^{k}\approx A$ for any $A\in\tco$.
Therefore by Theorem \ref{mainthonklus} (i), there is a natural isomorphism

\[
((A\otimes V_{0}^{k}\otimes D(A))_{\Gamma})^{d}\approx Hom_{\hgtg}(A,B), \]
for any $B\in\tco$.

 Recall that the  space
$((A\otimes V_{0}^{k}\otimes D(A))_{\Gamma})^{d}$ was defined by means of $\Gamma$, the
latter being defined by choosing three points,  $ \infty,1,0$,
see the end of \ref{Usingthespacesofcoinvariants}. The
choice of points was, of course, rather arbitrary. Keeping $ \infty, 0$
fixed and $A$, $D(B)$ attached to $\infty,0$ resp., we shall allow the third point to vary. We get then the family
of Lie algebras $\Gamma_{t}, t\in\nc^{\ast}$ and the family
of the  one-dimensional spaces (c.f. \ref{thesetup})
\[<A,V_{0}^{k},D(B)>_{t}:=
((A\times V_{0}^{k}\times D(B))_{\Gamma_{t}})^{d},\; t\in\nc^{\ast}. \]

These naturally arrange in a trivial line bundle over $\nc^{\ast}$,
the fiber being isomorphic to
\[<A,V_{0}^{k},D(B)>_{t}=(A\otimes D(B))_{\tgtg}=
Hom_{\hgtg}(A,B),\]
by the arguments using  Frobenius reciprocity as in  \ref{Usingthespacesofcoinvariants}. Pick a section of this bundle by choosing 
$\phi\in Hom_{\hgtg}(A,B)$.

 Hence we get a trilinear  functional (depending on $t\in\nc^{\ast}$)
\[\Phi_{t}^{\phi}\in <A,V_{0}^{k},D(B)>_{t}\subset (A\otimes V_{0}^{k}\otimes D(B))^{d}.\]

Reinterprete it as the linear map:

\begin{equation}
\label{fromcoinvtomaps_0}
\Phi_{t}^{\phi}(.): V_{0}^{k}\rightarrow (A\otimes  D(B))^{d},
\end{equation}

or, equivalently,
\begin{equation}
\label{fromcoinvtomaps}
\tilde{\Phi}_{t}^{\phi}(.): V_{0}^{k}\rightarrow Hom_{\nc}(A,  D(B)^{d}),\; t\in\nc^{\ast}.
\end{equation}

The latter map is an analogue of $N\rightarrow Hom_{\nc}(M,M\otimes N)$
 mentioned above.
To analyze its properties observe that there is an obvious embedding 
$B\rightarrow (D(B))^{d}$. It does not, of course,
 allow us to interprete $\tilde{\Phi}_{t}^{\phi}(v),\;v\in V_{0}^{k}$
as an element of $Hom_{\nc}(A,B)$ depending on $t$. But, as the following
lemma shows, Fourier coefficients of
 $\tilde{\Phi}_{t}^{\phi}(v),\;v\in V_{0}^{k}$ are actually elements
of $Hom_{\nc}(A,B)$. 
To formulate this lemma observe that there is a natural gradation
on $A$ and $B$ consistent with that of $\tgtg$; e.g.
$A=\oplus_{n\geq 0}A[n],\; dim A[n]<\infty$.

\begin{lemma}
\label{whycantfour}

Let $B$ be either $A$ or a quotient of $A$, $id: A\rightarrow B$ be the natural
projection. Then:

(i) $\Phi_{t}^{id}(vac)(x,y)=y(x),$ where $vac$ is understood as the generator
of $V_{0}^{k}$;

(ii) more generally, if $v\in V_{0}^{k}[n],\; x\in A[m],\;y\in D(B)[l]$,
then
\[\Phi_{t}^{id}(v)(x,y)\in \nc\cdot t^{-l+m-n}.\]
\end{lemma}

{\bf Proof.}  Given $g\in\gtg$, denote by $g_{n}\in\hgtg^{P}$ the
element $g\otimes (z-P)^{n}$ or $g\otimes z^{-n}$ if
$P=\infty$. (It should be clear from the
context which $P$ is meant.) Thus $g_{n}x=(g\otimes z^{-n})x$ if $x\in A$, the $A$
 being attached to $\infty$; similarly, 
  $g_{n}x=(g\otimes z^{n})x$ if $x\in D(B)$, the $D(B)$ being attached to $0$.

 (i) can be proved by an obvious induction on the
degree of $x$ and $y$ using the following formula (which follows from the definition
of  $((A\times V_{0}^{k}\times D(B))_{\Gamma})^{d}$ and the Laurent
expansions of $z^{-n}$ at $\infty$ and 0):

\[\Phi_{t}^{id}(vac)(g_{n}x,y)=-\Phi_{t}^{id}(vac)(x,g_{-n}y).\]

 To prove (ii) observe, first, that (i) is a particular case of (ii) when
$v=vac$. One then proceeds by induction on $n$ using the formula
 (which again follows from the definition
of  $((A\times V_{0}^{k}\times D(B))_{\Gamma})^{d}$
 and the Laurent
expansions of $(z-t)^{-n}$ at $\infty$ and 0):

\begin{multline}
(-1)^{n-1}(n-1)!\Phi_{t}^{id}(g_{-n}v)(x,y)\\
=(\frac{d}{dt})^{n-1}\{\sum_{i=1}^{\infty}t^{i-1}\Phi_{t}^{id}(v)(g_{i}x,y)
-\sum_{i=0}^{\infty}t^{-i-1}\Phi_{t}^{id}(v)(x,g_{i}y)\}.
\end{multline}

$\qed$

Observe that the spaces $A,B$ being graded, the space $Hom_{\nc}(A,D(B))$ is also.
Lemma \ref{whycantfour} means that although the map
$\tilde{\Phi}_{t}^{id}(.)$ from (\ref{fromcoinvtomaps}) cannot be interpreted
as an element of  $Hom_{\nc}(A,B)$, its Fourier components can  because they
are homogeneous. To compare with \cite{fren_lep_meur} introduce the following notation: for any $v\in V_{0}^{k}[n]$ set

\begin{equation}
Y(v,t)=\sum_{i\in\nz}v_{i} t^{-i-n},
\label{defofyvt}
\end{equation}
where
\begin{equation}
v_{i}:=\oint \tilde{\Phi}_{t}^{id}(v)t^{i+n-1}\, dt:
A[l]\rightarrow B[l+i],
\label{defofyvtfourcomp}
\end{equation}
for all $l\geq 0$, and call the generating functions $Y(v,t)$ {\em fields}.
 For example, it easily follows from the formulae above
that 
\begin{equation}
\label{current}
x(t):=Y(x_{-1}vac,t)=\sum_{i\in\nz}x_{i}t^{-i-1},
\end{equation}
producing the famous {\em current} $x(t)$. Another fact easily reconstructed
from the formulae above (especially from the proof of Lemma \ref{whycantfour})
is that 

\begin{equation}
\label{succcurrent}
(-1)^{n-1}(n-1)!Y(x_{-n}v,t)= :x(t)^{(n-1)}Y(v,t):,
\end{equation}

where we set
\[  :x(t)^{(n-1)}Y(v,z):\;=(x(z)^{(n-1)})_{-}Y(v,t) +
Y(v,t)(x(z)^{(n-1)})_{+},\]
$(x(z)^{(n-1)})_{\pm}$ being defined as usual (see e.g. \cite{frzhu}). 
 It follows that all fields are infinite
combinations of elements of $\hgtg$. 

The expressions  $Y(v,t)$ are not only formal generating functions.
In this notation Lemma \ref{whycantfour} can be rewritten as follows.

\begin{corollary}
\label{actionoffields3points}
Under the assumptions of Lemma \ref{whycantfour},
\[\Phi_{t}^{id}(v)(x,y) = y(Y(v,t)x).\]
\end{corollary}

\subsubsection{ Generalization   }
\label{twogeneralizations} The considerations of \ref{Vertexoperatorsandvertexoperatoralgebras}
 are easily generalized as follows. 
(We shall skip the proofs as they  essentially  repeat those in
\ref{Vertexoperatorsandvertexoperatoralgebras}.)

 Replace $V_{0}^{k}$ with $V_{\lambda}^{k}$ and pick $A,B\in\tco$
so that the space $<A,V_{\lambda}^{k},D(B)>_{t}\neq 0$. For any
$\phi\in <A,V_{\lambda}^{k},D(B)>_{t}$ we get a map

\begin{multline}
Y(.,t):\; V_{\lambda}^{k}\rightarrow Hom_{\nc}(A,B((t,t^{-1})))\\
V_{\lambda}^{k}\ni v\mapsto Y(v,t)=\sum_{i\in\nz}v_{i}t^{-i-\tilde{v}},\;
v_{i}\in  Hom_{\nc}(A,B).
\end{multline}

$Y(v,t),\; v\in V_{\lambda}^{k}$ is a generating function having all properties
its counterpart from \ref{Vertexoperatorsandvertexoperatoralgebras} 
with one notable exception. Consider the ``upper floor'' of $V_{\lambda}^{k}$:
$V_{\lambda}\subset V_{\lambda}^{k}$. The Fourier components
of the fields $Y(v,t),\; v\in V_{\lambda},\; \lambda\neq 0$ generate a $\hgtg$-submodule of 
$ Hom_{\nc}(A,B)$ isomorphic to the loop module
$L(V_{\lambda})=V_{\lambda}\otimes \nc[z,z^{-1}]$. Strange as it may seem to be,
if $\lambda=0$, then instead of $\nc[z,z^{-1}]$ this construction gives
simply $\nc$ -- this was explained above.

The embedding $L(V_{\lambda})\subset Hom_{\nc}(A,B)$ is called {\em a vertex operator}.
It is easy to see that all vertex operators are obtained via the described
construction.

\subsection{...and vertex operator algebras}
\label{and_vertex_operator_algebras}

We now recall that a vertex operator algebra (VOA) is defined to be a graded
vector space $\bigcup_{i\in\nz} V[i],\; dim V_{i}<\infty$ along with a map
\[Y(.,t): V\rightarrow End(V)((t,t^{-1})),\]
satisfying certain axioms among which we mention {\em associativity}
and {\em commutativity} axioms, see e.g. \cite{fren_lep_meur,frzhu}. Similarly
one defines the notion of a module (submodule) over a VOA. A VOA is a module over itself;
call an ideal of a VOA a submodule of a VOA as a module over itself. Observe
that it follows from the associativity axiom that the Fourier components of
fields $Y(v,t),\; v\in V$ close in a Lie algebra, $Lie(V)$. In this way,
an ideal of a VOA $V$ produces an ideal of $Lie (V)$ in the Lie algebra sense.
Not any ideal of $Lie(V)$ can be obtained in this way. Refer to such an ideal
{\em an ideal of }  $Lie(V)$ in the sense of {\em VOA}.

It follows from \cite{frzhu} that the constructions of 
\ref{Vertexoperatorsandvertexoperatoralgebras} give: 
$(V_{0}^{k}, Y(.,t))$ is a vertex
operator algebra and each $A\in\tco$ is a module over it. $Lie(V_{0}^{k})$
is habitually denoted $U(\hgtg)_{loc}$ and called {\em a local completion}
of $U(\hgtg)$, even though it is not an associative algebra! 
A moment's thought shows that the ideal lattice of
$U(\hgtg)_{loc}$ in the sense of
 VOA is isomorphic with the submodule lattice of
$V_{0}^{k}$
considered as a $\hgtg$-module.

\subsection{ }

Here we prove the following theorem -- one of the main results of
 this paper.

\begin{theorem}
\label{theoronannideals}
Let $k\in\nq_{>}$, $(\lambda,k),(0,k)\in P^{+}_{k}$ be regular.
Denote by $\Omega(V_{\lambda}^{k})$ the submodule lattice of
 $V_{\lambda}^{k}$,
and by $\Omega(U(\hgtg)_{loc})$ the ideal lattice of $U(\hgtg)_{loc}$ 
in the sense of VOA at the level $k$.
There is a lattice equivalence
\begin{equation}
\omega: \Omega(U(\hgtg)_{loc})\rightarrow \Omega(V_{\lambda}^{k}),\;\;
 \Omega(U(\hgtg)_{loc})\ni I\mapsto IV_{\lambda}^{k}.
\end{equation}
\end{theorem}

{\bf Proof.}

 First of all, by definition \ref{and_vertex_operator_algebras}  $\omega$ is
equivalently reinterpreted as a map of the submodule lattices
of the $\hgtg$-modules: $\omega:\; \Omega(V_{0}^{k})\rightarrow \Omega(V_{\lambda}^{k})$.
In what follows we shall make use of this reinterpretation.

  Consider the translation functor:
$T_{0}^{\lambda}$. If $N\subset V_{0}^{k}$ is a submodule,
then on the one hand we have
\[T_{0}^{\lambda}(V_{0}^{k})=V_{\lambda}^{k}\tens V_{0}^{k} (=V_{\lambda}^{k}),\]
and therefore
\[T_{0}^{\lambda}(V_{0}^{k}/N)=V_{\lambda}^{k}\tens ( V_{0}^{k}/N).\]

By Theorem \ref{mainthonklus} and Corollary \ref{actionoffields3points}, 
\begin{multline}
Hom_{\hgtg}(T_{0}^{\lambda}(V_{0}^{k}/N),?)= <V_{0}^{k}/N,V_{\lambda}^{k},D(?)>_{t}\\
= Hom_{\hgtg}(V_{\lambda}^{k}/\omega(N),?).
\end{multline}

On the other hand, by Proposition \ref{conseqreltochar} (i)
 \[T_{0}^{\lambda}(V_{0}^{k}/N)= V_{\lambda}/T_{0}^{\lambda}(N).\]

We conclude immediately that $\omega(N)= T_{0}^{\lambda}(N)$.
It remains to recollect that $T_{0}^{\lambda}$ is an isomorphism of
the submodule lattices by Proposition \ref{conseqreltochar} (v).
$\qed$

\bigskip

An application of this result to annihilating ideals of
{\em admissible representations} is as follows. Recall that if
$k+h^{\vee}\in\nq_{>}$, $(\lambda,k)\in P^{+}_{k}$ is regular, then
$L_{\lambda}^{k}$ is called admissible \cite{kac_wak}.
$L_{\lambda}^{k}$ is an irreducible quotient of
$V_{\lambda}^{k}$ be a submodule $N_{\lambda}^{k}$
 generated by one singular vector,
see also \cite{kac_wak}.  By Theorem \ref{theoronannideals},
$\omega(N_{0}^{k})=N_{\lambda}^{k}$. We get 

\begin{corollary}
\label{annidadmrepr}
The annihilating ideal of an admissible representation  equals
$Lie(N_{0}^{k})$; in particular, it is generated (in the sense of VOA) by one
singular vector of $V_{0}^{k}$.
\end{corollary}

{\em Remarks.}
(i) Theorem \ref{theoronannideals} reduces the problem of classifying
annihilating ideals  to the easier problem of classifying 
submodules of $V_{0}^{k}$. What can be said about the latter?
It has been known for a while that in the simplest case of
$\hat{sl}_{2}$ $V_{0}^{k}$ contains a unique proper non-trivial
submodule. In general, multiplicities in composition series
are described by Kazhdan-Lusztig polynomials. In the recent work
\cite{lusztig} Lusztig exhibits an example when infinitely many
multiplicities are non-zero and thus there are infinitely many
different submodules. This result makes one believe that this is 
" usually" the case.

(ii) In the case $\gtg=\gtsl_{2}$, Corollary \ref{annidadmrepr}
follows from the more general results of \cite{feimal1}, see also
\cite{feimal2}.

(iii) If the Feigin-Frenkel conjecture on the singular support of
$L_{0}^{k}$ (theorem in the $sl_{2}$-case, see \cite{feimal2})
were correct, then Corollary \ref{annidadmrepr}
 would imply its validity for any admissible representation from
$\tco$ and thus would give a new example of rational conformal field
theory.

(iv) Another way to think of Corollary \ref{annidadmrepr} is
that $L_{0}^{k}$ is a VOA and $L_{\lambda}^{k}$ is a module over
it; in the $sl_{2}$-case, this point of view is adopted in
 \cite{adam,donglimas}.

\section{ tilting functors -- the negative level case }
\label{whcafm}
From now on $k+h^{\vee}$
 is  a negative rational number. The structure of $\co_{k}$
can be described as follows. Consider an antidominant 
weight $(\lambda,k)$ with integral $\lambda$.
 Denote by $\co_{k}^{\lambda}$ the
full subcategory consisting of 
modules containing only irreducibles with highest weights
lying in the $W_{k}$-orbit of an antidominant weight $(\lambda,k)$. 
By \cite{deodgabbkac,rcw},     there is the decomposition
$\co_{k}=\oplus_{\lambda}\co_{k}^{\lambda}$ (c.f. \ref{preliminaries}).

\begin{sloppypar}
It follows from \cite{kac_kazhd} that any module
 from $\co_{k}^{\lambda}$ has a finite length. There
arises the Grothendieck ring,  each of   the
following   sets $\{ch(M_{w\cdot\lambda}^{k}),\; w\in W_{k}\}$,  
 $\{ch(L_{w\cdot\lambda}^{k}),\; w\in W_{k}\}$ being a basis of it.
\end{sloppypar}

All this carries over to the case of $\tco$ 
by replacing $M_{\lambda}^{k}$ with 
$V_{\lambda}^{k}$ and making sure that
 in the latter case $\lambda$ is dominant integral.

\subsection{Tilting Modules}
\label{clasoftitltsect}

{\bf Definition} 
A module $W$ from $\co_{k}$ 
($\tco$ resp.) is called {\em tilting }if both
$W$ and $D(W)$ possess a filtration by Verma (Weyl resp.) modules. 

\bigskip

For instance,
a Verma module
with an antidominant highest weight is
  a tilting module in $\co_{k}$, for it is
irreducible and, therefore, isomorphic to its dual. Likewise, a Weyl module
with an antidominant highest weight is  a tilting module in $\tco$.

\begin{proposition}
\label{classtiltmod}

(i) Any tilting module from $\tco$
 is a direct sum of indecomposable tilting modules.
For any $(\mu,k)$, $\mu$ being dominant integral,
 there is a uniquely determined
indecomposable tilting module $\tilde{W}_{\mu,k}$ such that

\[ch(\tilde{W}_{\mu,k})=ch(V_{\mu}^{k})+ \sum_{\nu<_{k}\mu}c_{\nu} 
ch(V_{\nu}^{k}).\]

(For definition of $<_{k}$ see \ref{defofordandlem}.) 
The map $\mu\mapsto \tilde{W}_{\mu,k}$ establishes a bijection between the set of weights 
satisfying the above mentioned condition and
the set of isomorphism classes of indecomposable tilting modules in $\tco$.

(ii) Likewise, any tilting module
 from $\co_{k}$ is a direct sum of indecomposable
 tilting modules.
For any $(\mu,k)$, there is a uniquely
 determined indecomposable tilting module $W_{\mu,k}$ such that

\[ch(W_{\mu,k})=ch(M_{\mu}^{k})+ 
\sum_{\nu<_{k}\mu}c_{\nu} ch(M_{\nu}^{k}).\]

The map $(\mu,k)\mapsto W_{\mu,k}$ establishes
 a bijection between the set of weights and
the set of isomorphism classes of indecomposable tilting modules in $\co_{k}$.
\end{proposition}

{\bf Proof} (i) is proved in \cite{kazh_luszt_1}
 using general results of Ringel \cite{ringel}; 
(ii) can be proved in   the same way. A lucid exposition suited
for the $\co$-category case can also be found in \cite{and}. $\qed$

\bigskip

\begin{corollary}
\label{basistilt}
The set of characters of tilting modules
 in $\co_{k}$ ($\tco$ resp.) is a basis
of the Grothendieck ring of $\co_{k}$ ($\tco$ resp.) 
\end{corollary}

{\bf Proof} -- obvious.

\bigskip

Kazhdan and Lusztig prove that the category $\tco$ is rigid.
Some of the consequences of this fact are collected in the following 

\begin{lemma}
\label{consofrig_6}

(i) For any $A\in\tco$, the functor
 $A\tens ?: \co_{k}\rightarrow \co_{k}$ is
exact.

(ii) For any $A\in\tco$, $B\in \co_{k}$,
 one has: $D(A\tens B)=D(A)\tens D(B)$.

(iii) If $A\in\tco$, $B\in \co_{k}$ are both tilting, then $A\tens B$ is also.
\end{lemma}

{\bf Proof} (i) is proved in the
 same way Corollary \ref{adjointnessproperty} was proved.

(ii) is a general fact about monoidal
 categories, see \cite{kazh_luszt_1} III, Proposition A.1.

(iii) is an immediate consequence of
 (ii) and Lemma \ref{multinfiltr} which is valid along
with its proof for all values of central charge. $\qed$

\bigskip

\subsection {Tilting Functors}
\subsubsection{{\em Modules depending on a parameter}}
\label{newcategories} 
Let $t$ be an indeterminate and $R(t)$ the ring of
rational functions having poles only on a positive real ray. Along with
$\co_{k}$, introduce the category $\co_{t}$ defined in the same way
 except that $\nc$ as a ground field is replaced with $R(t)$, that is,
$\hgtg$
is regarded as an $R(t)$-algebra and objects of $\co_{t}$ are required
to be free $R(t)$-modules. There is the specialization functor
\[sp_{k}:\co_{t}\rightarrow\co_{k},\; M\mapsto
M(k)\stackrel{\mbox{def}}{=}M/(t-k)M.\]

Let $\ca_{t}\subset\co_{t}$ be the full subcategory of modules $V$
such that $V(k)$ has a Weyl filtration. Define
\[^{fam}\co_{k}\stackrel{\mbox{def}}{=}sp_{k}(\ca_{t}).\]

In other words, $^{fam}\co_{k}$ consists of modules having a Verma
filtration and allowing inclusion in a family analytically depending
on $k$.

One defines in a similar way $^{fam}\co^{\lambda}_{k}\subset
\co^{\lambda}_{k}$ and the Weyl versions 
$^{fam}\tco$, $^{fam}\tco^{\lambda}\subset ^{fam}\tco$.

\bigskip
{\em Examples}

(i) $V_{\lambda}^{k}\in^{fam}\tco$, $M_{\lambda}^{k}\in^{fam}
\co_{k}$.

(ii) $V_{\lambda}^{k}\tens V_{\mu}^{k}\in^{fam}\tco$,
 $V_{\mu}^{k}\tens M_{\lambda}^{k}\in^{fam}
\co_{k}$.

(iii) Direct summands of any of the modules from (ii) belong
to the corresponding categories.

\bigskip

\subsubsection{{\em Main Results}}

From now on we shall have two  fixed  regular antidominant weights
 $(\lambda_{r},k)$ and
$(\lambda_{l},k)$ so that $\lambda_{r}-\lambda_{l}$ is integral.

{\bf Definition} Tilting functor is a direct summand of the functor
\[p_{\lambda_{l}}\circ(W\tens ?): ^{fam}\co_{k}^{\lambda_{r}}\rightarrow
^{fam} \co_{k}^{\lambda_{l}},\]
where $W\in \tco$ is tilting.

\bigskip

The following is one of our main results.

\begin{theorem}
\label{class_titl_fun}

(i) Each tilting functor is a direct sum of indecomposable ones.

(ii) There is a 1-1 correspondence between tilting functors and elements
of $W_{k}$. The functor $\Phi_{w}\in Funct(^{fam}\co_{k}^{\lambda_{r}},
^{fam}\co_{k}^{\lambda_{l}})$
attached to $w\in W_{k}$ is uniquely determined by the condition that
\[\Phi_{w}(M_{\lambda_{r}}^{k})= W_{w\cdot\lambda_{l},k}.\]
\end{theorem}

\begin{sloppypar}
We shall prove this theorem in  \ref{proofs}. 
One consequence of this theorem is that the category
of tilting functors, as a subcategory of
 $Funct(\co_{k}^{\lambda_{r}},\co_{k}^{\lambda_{l}})$,
 is semi-simple and isomorphic
to the category of tilting modules in $\co_{k}^{\lambda_{l}}$. 
Denote this category by
$Tilt(\lambda_{l},\lambda_{r})$.
\end{sloppypar}

\begin{sloppypar}
Suppose now that we are given three antidominant 
highest weights, $ (\lambda_{1},k),
(\lambda_{2},k), (\lambda_{3},k)$, and two tilting functors: 
$\Phi\in Funct(^{fam}\co_{k}^{\lambda_{1}}, 
^{fam}\co_{k}^{\lambda_{2}} )$ and 
$\Psi\in Funct(^{fam}\co_{k}^{\lambda_{2}},^{fam}\co_{k}^{\lambda_{3}}) $.
 Then $\Psi\circ\Phi$
is also  a tilting functor.
 This follows from Lemma \ref{consofrig_6}(iii) and the
 associativity morphism:
 $W_{2}\tens (W_{1}\tens B)=(W_{2}\tens W_{1})\tens B$.
If in addition $\lambda_{l}=\lambda_{r}=\lambda$,
 then the category $ Tilt(\lambda) \stackrel
{def}{=}Tilt(\lambda_{r},\lambda_{l})$ is closed under $\circ$ and the pair 
$(Tilt(\lambda), \circ)$ becomes a semi-simple
 monoidal category. The following is
another main result concerning tilting functors.
\end{sloppypar}

\begin{theorem}
\label{grothringoftilt}
The Grothendieck ring of
 $(Tilt(\lambda), \circ)$ is isomorphic to the group algebra
of $W_{k}$.
\end{theorem}

Proof of this theorem is to be found in the next section.

\subsection {Further theorems on tilting functors and proofs}
\label{proofs}

\subsubsection{ }
\label{charintermofimage}
Analogously to the semi-simple case
 (see Theorems \ref{classprojbg} and \ref{charprojmv}),
the key to the proof of 
Theorem \ref{class_titl_fun} is the following result.

\begin{theorem}
\label{titlfunvsitsvalue}
Let $\Phi,\Psi\in 
Funct(\co_{k}^{\lambda_{r}},\co_{k}^{\lambda_{l}})$ be tilting functors.
The natural map 
\[Mor (\Phi,\Psi)\rightarrow Hom_{\hgtg}(\Phi(M_{\lambda_{r}}^{k}),\Psi(M_{\lambda_{r}}^{k}))\]
is a vector space  isomorphism.
\end{theorem}

{\bf Proof} It is obviously enough 
to consider the case when $\Phi=A_{1}\tens ?$,
$\Psi=A_{2}\tens ?$, 
$A_{1},A_{2}$ being tilting modules. As $\tco$ is rigid, we have
\[Hom_{\hgtg}(A_{1}\tens ?, A_{2}\tens ?)
=Hom_{\hgtg}(D(A_{2})\tens A_{1}\tens?, ?).\]
Therefore it is enough to prove that the natural map
\begin{equation}
\label{reduction}
 Mor (A\tens ?, id) \rightarrow 
Hom_{\hgtg}(   A\tens M_{\lambda_{r}}^{k},M_{\lambda_{r}}^{k})
\end{equation}
is an isomorphism for any tilting $A$.

\bigskip

{\bf Injectivity of (\ref{reduction})} Let $\phi\in Mor ( A\tens ?, id) $
be such that its value $\phi(id)$ on 
$id\in End ( M_{\lambda_{r}}^{k} )$ is zero. 
 Consider the exact
sequence
\[0\rightarrow M_{\lambda_{r}}^{k}
 \rightarrow M_{u\cdot\lambda_{r}}^{k} \rightarrow
 M_{u\cdot\lambda_{r}}^{k} /M_{\lambda_{r}}^{k}\rightarrow 0 .\]
By Lemma \ref{consofrig_6} (i), the sequence

\[ 0\rightarrow A\tens  M_{\lambda_{r}}^{k} \rightarrow 
  A\tens M_{u\cdot\lambda_{r}}^{k} \rightarrow
   A\tens (M_{u\cdot\lambda_{r}}^{k} /M_{\lambda_{r}}^{k})\rightarrow 0 \]
is also exact. Therefore, 
\[\phi(id)\in Hom_{\hgtg}(  A\tens M_{u\cdot\lambda_{r}}^{k},
 M_{u\cdot\lambda_{r}}^{k})\]
factors through to the map
\[  A\tens (M_{u\cdot\lambda_{r}}^{k} /M_{\lambda_{r}}^{k}) \rightarrow
M_{u\cdot\lambda_{r}}^{k}.\]
To show that this map can only be zero we make use of the notion of
singular support.

Recall that under certain technical assumptions (which are satisfied
in the $\co$-category case,see e.g.
\cite{feimal2} for the necessary definitions) one defines the singular
support of a Lie algebra module to be the zero set of the annihilating
ideal of the corresponding graded object; thus singular support is
a conical subset of the dual to the Lie algebra in question. For example,
the singular support of a Verma module is all functionals vanishing
on the Borel subalgebra. If a module has a finite filtration, then
its singular support is the union of the singular supports of the
successive quotients. Thus any module with a Verma filtration,
$A\tens M_{u\cdot\lambda_{r}}^{k}$ for example, has
the same support as a Verma module.

A quotient of a Verma module by a proper submodule,
 however, has a smaller singular support;
it is obtained  by imposing additional equations, those coming from
the symbols of the submodule. By the same token,  the singular support of 
$  A\tens (M_{u\cdot\lambda_{r}}^{k} /M_{\lambda_{r}}^{k}) $ is strictly  less
than that of a Verma module. Therefore, the latter module may not
contain a Verma module as a subquotient and the map
\[  A\tens (M_{u\cdot\lambda_{r}}^{k} /M_{\lambda_{r}}^{k}) \rightarrow
M_{u\cdot\lambda_{r}}^{k}\]
can only be zero.
Thus $\phi(id)$ is zero on
any Verma module and, by exactness,
 is also on any module from $\co_{k}$. Injectivity has been
proven.

\bigskip

{\bf Surjectivity of (\ref{reduction})} What remains to be done
is to show that any element of
 $Hom_{\hgtg}(A\tens M_{\lambda_{r}}^{k},M_{\lambda_{r}}^{k})$
naturally determines an element of 
$Hom_{\hgtg}(A\tens B,B) $ for any $B\in^{fam}\co_{k}^{\lambda_{r}}$
 and any tilting $A\in\tco$. 
We begin with calculating the space 
 $Hom_{\hgtg}(A\tens M_{\lambda_{r}}^{k},M_{\lambda_{r}}^{k})$.
Recall that $A$ has a Weyl filtration with quotients $V_{\mu_{i}}^{k}$,
$i=1,...,n$.
\begin{lemma}
\label{vermawithantdom}
\[Hom_{\hgtg}(A\tens M_{\lambda_{r}}^{k},M_{\lambda_{r}}^{k})=
\oplus_{i=1}^{n}V_{\mu_{i}}[0],\]
where $V_{\mu_{i}}[0]$ stands for the zero weight subspace of the
$\gtg$-module $V_{\mu_{i}}$. 
\end{lemma}

{\bf Proof} 
Let, first, $A=V_{\mu}^{k}$.  Then by Theorem \ref{mainthonklus}
we have
\begin{multline}
Hom_{\hgtg}(V_{\mu}^{k}\tens M_{\lambda_{r}}^{k},M_{\lambda_{r}}^{k})
=[(V_{\mu}^{k}\otimes M_{\lambda_{r}}^{k}\otimes D(M_{\lambda_{r}}^{k}))_
{\Gamma}]^{d}\\
=[(V_{\mu}^{k}\otimes M_{\lambda_{r}}^{k}\otimes M_{\lambda_{r}}^{k})_
{\Gamma}]^{d}=V_{\mu}[0],
\end{multline}
where the second equality follows from the isomorphism 
$M_{\lambda_{r}}^{k}=D(M_{\lambda_{r}}^{k}$, while the third one
follows from the Frobenius reciprocity (c.f. the derivation of 
(\ref{morphofweylintodual}) in \ref{Usingthespacesofcoinvariants}.

To treat the case when $A$ has a Weyl filtration with more than one
term, observe that $A\tens M_{\lambda_{r}}^{k}$ has a Verma filtration.
As in the proof  of
Proposition 20.1 (6) in \cite{kazh_luszt_1} III sect.20, one shows that
\[Ext^{1}(B,D(M_{\lambda_{r}}^{k}))=0\]
for any $B$ carrying a Verma filtration. The long cohomology sequence
implies then that the space
\[Hom_{\hgtg}(B,D(M_{\lambda_{r}}^{k}))\]
behaves as if $B$ were a direct sum of Verma modules. Lemma is proved.
$\qed$

\bigskip
Essential for the proof of Lemma \ref{vermawithantdom}
were the following properties of the modules in question:
$A\tens M_{\lambda_{r}}^{k}$ is a free $\hat{\gtnm}$-module
(this is equivalent to carrying 
a Verma filtration), and $D(M_{\lambda_{r}}
^{k})$ is a co-induced module, that is, dual to an induced one.
Therefore, the same arguments give the following more general
result.

\begin{lemma}
\label{calcmorphindcoind}
Let $\gtg=\gtnm\oplus\gth\oplus\gtnp$ and
$\hgtg=\hat{\gtnm}\oplus\hat{\gth}\oplus\hat{\gtnp}$
be corresponding triangular decompositions. Further,
let $C\in\co_{k}$ be freely generated by $\hat{\gtnm}$
from a finite dimensional $\gth$-space $V$ and
\[B=Ind_{\hat{\gth}\oplus\hat{\gtnp}}^{\hgtg}W,\; dim W<\infty.\]
Then
\[Hom_{\hgtg}(C,D(B))=Hom_{\hat{\gth}\oplus\hat{\gtnp}}
(C,D(W))=Hom_{\gth}(V,D(W)).\]
\end{lemma}

What does not allow us to extend Lemma\ref{vermawithantdom}
to all $B\in\co^{k}$ is that in general $B\neq D(B)$. Let,
however, $k$ be generic, that is, irrational. Then $\tco$
is semi-simple and $A=\oplus_{i}V_{\mu_{i}}^{k}$; $\co_{k}$
is equivalent to $\co_{\gtg}$, the $\co$-category of $\gtg$-modules;
under this equivalence $A$ goes to$\oplus_{i}V_{\mu_{i}} $
and $\tens$ is transformed $\otimes$, see \cite{fink}. Therefore
in this case, the surjectivity of (\ref{reduction}) becomes the
corresponding statement of the  "semi-simple theory"  -- Theorem
3.5 of
\cite{bernst_gelf}.

Now the strategy   becomes obvious: include $A$ and $B$ in a family
of modules depending on $k\in\nc$ and prove that all morphisms
existing generically admit continuation to our particular value of $k$;
it is at this point we use the fact that $B\in^{fam}\co_{k}$,
$A\in ^{fam}\tco$.  It suffices to prove the following lemma.

\begin{lemma}
\label{howtocontinue}
Let ${\bf A}\in\tilde{\co}_{t},\; {\bf B}\in\co_{t}$ be such that
${\bf A}(k)=A,\; {\bf B}(k)=B$. Then

(i) $Hom_{\co_{t}}({\bf A}\tens {\bf B},{\bf B})$ is a free
$R(t)$-module;

(ii) the natural map
\[Hom_{\co_{t}}({\bf A}\tens {\bf B},{\bf B})(k)\rightarrow
Hom_{\co_{k}}( A\tens  B, B)\]
is an embedding.
\end{lemma}

{\bf Proof} follows the lines of the proof of Lemma 22.8 in 
\cite{kazh_luszt_1}.

(i) If ${\bf B}=D({\bf B'})$ for some induced representation 
  ${\bf B'}$, then (i) follows from Lemma \ref{calcmorphindcoind}.
In general, there is an induced representation ${\bf B'}$ such that
${\bf B'}$ maps onto $D({\bf B})$ and therefore ${\bf B}$ embeds
in $D({\bf B'})$. Hence  $Hom_{\co_{t}}({\bf A}\tens {\bf B},{\bf B})$
is a submodule of a free $R(t)$-module and thus is also free.

(ii) Let $h\in Hom_{\co_{t}}({\bf A}\tens {\bf B},{\bf B})(k)$
be such that the image of $h(k)$ in $Hom_{\co_{k}}( A\tens  B, B)$ is zero.
It means that $h$ evaluated on any element of ${\bf A}\tens {\bf B}$
belongs to $(t-k){\bf B}$. Therefore 
$h\in (t-k)Hom_{\co_{t}}({\bf A}\tens {\bf B},{\bf B})$ and hence
$h(k)=0$. $\qed$ 

Theorem \ref{titlfunvsitsvalue} has been proved. $\qed$

\subsubsection{Proof of Theorem \ref{class_titl_fun} }
By Theorem \ref{titlfunvsitsvalue}, tilting functors
$\Phi_{1},\Phi_{2}$ are isomorphic
if and only if $\Phi_{1}(M_{\lambda_{r}}^{k})\approx
\Phi_{2}(M_{\lambda_{r}}^{k})$ and direct summands
of any $W\tens ?$ are in 1-1 correspondence with direct
summands of $W\tens M_{\lambda_{r}}^{k}$. By Lemma \ref{consofrig_6},
$W\tens M_{\lambda_{r}}^{k}$ is tilting. Proposition \ref{classtiltmod}
implies that for any indecomposable tilting functor $\Phi$ there is
$w\in W_{k}$ so that
\[ \Phi(M_{\lambda_{r}})=W_{w\cdot\lambda_{l},k}.\]

It remains to prove that for any $w\in W_{k}$ there is a $\Phi$
so that $\Phi(M_{\lambda_{r}})= W_{w\cdot\lambda_{l},k}$. Pick a
dominant integral
$\mu\in\gth^{\ast}$ such that $w\cdot\lambda_{l}-\lambda_{r}$
is an extremal weight of $V_{\mu}$. Consider tilting module
$\tilde{W}_{\mu,k}\in\tco$.  
It follows from Lemma \ref{multinfiltr}  that the composition series of
 $\tilde{W}_{\mu,k}\tens M_{\lambda_{r}}^{k}$
contains $M_{w\cdot\lambda_{l}}^{k}$ and 
 by  Lemma \ref{longowbg} $w\cdot\lambda_{l}$ is maximal
among $\nu$ such that $M_{\nu}^{k}$ appears in the composition
series of  $\tilde{W}_{\mu,k}\tens M_{\lambda_{r}}^{k}$.
Proposition \ref{classtiltmod} immediately gives then that
$\tilde{W}_{\mu,k}\tens M_{\lambda_{r}}^{k}$
contains a unique direct summand isomorphic to
$W_{w\cdot\lambda_{l},k}$. $\qed$

\subsubsection{ }
 Being exact by Lemma \ref{consofrig_6},
each tilting functor induces a  
homomorphism
of the Grothendieck rings. On the other hand, the
 Grothendieck rings carry 
a natural action
of $W_{k}$.
 
\begin{theorem}
\label{wafflinearity}
The homomorphism of the Grothendieck rings induced 
by a tilting functor is $W_{k}$-linear.
\end{theorem}

{\bf Proof} For any exact functor $\Phi$ denote by $\Phi^{k}$ the 
induced homomorphism of Grothendieck rings
 ( or $K$-rings, hence the
notation). If $W\in\tco$, then $(W\tens ?)^{K}$ is 
$W_{k}$-linear as it equals the operator of multiplication
by the formal character of the corresponding finite dimensional
$\gtg$-module, see Lemma \ref{multinfiltr}.
Therefore it suffices to prove that given a tilting functor $\Phi$,
$\Phi^{K}$ equals a linear combination of functors $W\tens ?$
with tilting $W\in\tco$.

Let $\Phi$ be such that 
\[\Phi(M_{\lambda_{r}}^{k})=W_{w\cdot\lambda_{l},k},\; w\in W_{k}.\]
Pick a
dominant integral
$\mu\in\gth^{\ast}$ such that $w\cdot\lambda_{l}-\lambda_{r}$
is an extremal weight of $V_{\mu}$. Consider $\tilde{W}_{\mu,k}\tens ?$.
As we saw when proving Theorem\ref{class_titl_fun},
its direct sum decomposotion is of the form
\[\tilde{W}_{\mu,k}\tens ?=\Phi_{w}\oplus \oplus_{v<w}c_{v}\Phi_{v}.\]
Hence
\[(\tilde{W}_{\mu,k}\tens?)^{K}
=(\Phi_{w})^{K}+\sum_{v<w}c_{v}(\Phi_{v})^{K}.\]
 
The obvious induction on the length of $w$ allows us to  "solve" 
the latter equality for $(\Phi_{w})^{K}$ in terms of 
$(\tilde{W}_{v\cdot\mu,k}\tens?)^{K}$, $v\in W_{k}$.  $\qed$

\subsubsection{Proof of Theorem \ref{grothringoftilt} }
 Given tilting functor $\Phi_{w},\; w\in W_{k}$ (see 
Theorem \ref{class_titl_fun} ),
set $q_{w}=ch(\Phi_{w}(M_{\lambda}^{k}))$ and regard
 it as an element of the group
algebra of $W_{k}$.  By Theorem \ref{wafflinearity}, 
\[ch(\Phi_{w}(  M_{u\cdot\lambda}^{k}) )=u q_{w},\;u\in W_{k}.\]
As the set  $\{ch(M_{u\cdot\lambda}^{k}),\; u\in W_{k}\}$ 
is a basis of the Grothendieck
ring of $\co_{k}^{\lambda}$, we get that 
the action of $\Phi_{w}$ is  multilication by
$q_{w}$. To complete the proof it is enough to remark that
by Corollary \ref{basistilt}, the set $\{q_{w},\; w\in W_{k}\}$ 
is a  basis
of the group algebra of $W_{k}$. $\qed$


\begin{thebibliography}{99}

\bibitem{adam} Adamovic D., Milas A. MRL {\bf 2} (1995) 563-575

\bibitem{and} Andersen H., Paradowski J. CMP {\bf 169} (1995) 563-588
 
\bibitem{bernst_gelf} Bernstein J.N., Gelfand S.I.
Compositio Mathematica {\bf 41}, 2, (1980) 245-285

 
\bibitem{deodgabbkac} Deodhar V.V., Gabber O., Kac V.G., Adv.in Math. {\bf 45}
(1982) 92-116

 

\bibitem{drinf_1} Drinfeld V.G.   Algebra Anal. {\bf 2} (1990) 149-181

\bibitem{donglimas} Dong C., Li H., Mason G., Vertex operator algebras
associated to admissible representations of $\widehat{sl}_{2}$, 
q-algebra/9509026

\bibitem{duflo} Duflo M. Lect. Notes in Math. {\bf 497} (1975) 26-88

\bibitem{enr} Enright T. Ann. of Math. {\bf 110} (1979) 1-82


 
 \bibitem{feimal1} Feigin B., Malikov F. Lett.in Math.Phys.
{\bf 31} (1994) 315-325

 \bibitem{feimal2} Feigin B., Malikov F.
 Cont.Math. {\bf 202}  ``Operads: Proceedings of Renaissance Conferences''
(ed. by Loday, Stasheff and Voronov)
 

 \bibitem{fink} Finkelberg M. Fusion categories, Ph.D. thesis, Harvard university,1993

\bibitem{fkw} Frenkel E., Kac V., Wakimoto M. Comm.Math.Phys. {\bf 147}(1992) 295-328



 
\bibitem{fren_lep_meur} Frenkel I.B., Lepowsky J., Meurman A.
 {\em Vertex Operator
 Algebras and the Monster}, Academic Press, Inc 1988
 
\bibitem{frzhu} Frenkel I.B., Zhu Y., Duke Math. Journal {\bf 66} (1992) 123-168

\bibitem{jantz0} Jantzen J.C. Lect. Notes in Math. (1970)
 
\bibitem{jantz} Jantzen J.C. Representations of Algebraic Groups, Pure
and Applied Mathematics {\bf 131}, Academic Press (1987)
\bibitem{joseph} Joseph A. J. London Math. Soc. (2) {\bf 20} (1979),
193-204
 
\bibitem{kac_kazhd} Kac V.G., Kazhdan D.A., Adv.in Math.
{\bf 34} (1979) 97-108


\bibitem{kac_wak} Kac V.G., Wakimoto M. Proc. Nat'l Acad. Sci. USA
{\bf 1988} 4956

\bibitem{kazh_luszt_0} Kazhdan D., Lusztig G. Duke Math.J. {\bf 62} 21-29
\bibitem{kazh_luszt} Kazhdan D., Lusztig G.   JAMS {\bf 6} no. 4 (1993)

\bibitem{kazh_luszt_1} Kazhdan D., Lusztig G.   JAMS {\bf 6} no. 5 (1993)
\bibitem{lusztig} Lusztig G. Electronic Representation Theory {\bf 1}
(1997) 25-30

\bibitem{mirc} Mircovi\v{c} I., unpublished

\bibitem{ringel} Ringel C.M. Math.Z. {\bf 208} (1991), 209-223 
 

 
\bibitem{rcw} Rocha A., Wallach N. Math.Z. {\bf 180} (1982) 151-177

 
 
 \bibitem{vogan} Vogan D. Duke Math.J. {\bf 46} (1979)

\bibitem{zuck} Zuckerman G. Ann. of Math. {\bf 106} (1977) 295-308

 \end{thebibliography}
\end{document}